\documentclass{elsart3-1-ppt}

\usepackage{amssymb,amsmath,latexsym,setspace}
\usepackage[english,francais]{babel}

\newtheorem{e-proposition}[theorem]{Proposition}

\newtheorem{e-definition}[theorem]{Definition\rm}

\renewcommand{\theequation}{\arabic{equation}}
\newcommand{\M}{\mathfrak M}
\newcommand{\iM}[1]{\int_{\M}{#1}\,d\kern1pt v_g}

\newcommand{\Lap}{\Delta_g}
\newcommand{\Ric}{\mathfrak R}

\def\og{\leavevmode\raise.3ex\hbox{$\scriptscriptstyle\langle\!\langle$~}}
\def\fg{\leavevmode\raise.3ex\hbox{~$\!\scriptscriptstyle\,\rangle\!\rangle$}}

\newcommand{\R}{{\mathbb R}}
\renewcommand{\S}{{\mathbb S}}
\newcommand{\N}{{\mathbb N}}
\newcommand{\be}[1]{\begin{equation}\label{#1}}
\newcommand{\ee}{\end{equation}}
\renewcommand{\(}{\left(}
\renewcommand{\)}{\right)}

\newcommand{\seq}[2]{({#1}_{#2})_{#2\in\N}}

\usepackage{color}

\renewcommand{\S}{\mathbb{S}}

\newcommand{\nrm}[2]{\|{#1}\|_{\L^{#2}(\S^d)}}
\newcommand{\nrmM}[2]{\|{#1}\|_{\L^{#2}(\M)}}
\newcommand{\nrmRd}[2]{\|{#1}\|_{\L^{#2}(\R^d)}}

\renewcommand{\H}{\mathrm H}
\renewcommand{\L}{\mathrm L}

\newcommand{\muM}{\mu} 
\newcommand{\alphaM}{\alpha} 
\newcommand{\kappaS}{\kappa} 
\newcommand{\kappaM}{\kappa} 

\newcommand{\finprf}{\unskip\null\hfill$\;\square$\vskip 0.3cm}
\newenvironment{proof}{\par\noindent{\emph{Proof. }}}{\finprf}

\journal{the Acad\'emie des sciences}
\begin{document}

\centerline{}
\begin{frontmatter}

\selectlanguage{english}


\title{Spectral properties of Schr\"odinger operators on compact manifolds: rigidity, flows, interpolation and spectral estimates}

\selectlanguage{english}
\author[Ceremade]{Jean Dolbeault}
\ead{dolbeaul@ceremade.dauphine.fr}
\and
\author[Ceremade]{Maria J.~Esteban}
\ead{esteban@ceremade.dauphine.fr}
\and
\author[ICL]{Ari Laptev}
\ead{a.laptev@imperial.ac.uk}
\and
\author[GIT]{Michael Loss}
\ead{loss@math.gatech.edu}

\address[Ceremade]{Ceremade (UMR CNRS no. 7534), Universit\'e Paris-Dauphine, Place de Lattre de Tassigny, 75775 Paris 16, France}
\address[ICL]{Department of Mathematics, Imperial College London, Huxley Building, 180 Queen's Gate, SW7 2AZ, UK}
\address[GIT]{School of Mathematics, Skiles Building, Georgia Institute of Technology, Atlanta GA 30332-0160, USA}

\begin{abstract}\selectlanguage{english}
This note is devoted to optimal spectral estimates for Schr\"odinger operators on compact connected Riemannian manifolds without boundary. These estimates are based on the use of appropriate interpolation inequalities and on some recent rigidity results for nonlinear elliptic equations on those manifolds.
\vskip 0.5\baselineskip\noindent

\selectlanguage{francais}
\vskip 0.5\baselineskip \noindent {\bf Propri\'et\'es spectrales d'op\'erateurs de Schr\"odinger sur des vari\'et\'es compactes : rigidit\'e, flots, interpolation et estimations spectrales.}

\vskip 0.5\baselineskip\noindent{\bf R\'esum\'e.}
Cette note est consacr\'ee \`a des estimations spectrales optimales pour des op\'erateurs de Schr\"odinger sur des vari\'et\'es Riemaniennes compactes et simplement connexes, sans bord. Ces estimations sont bas\'ees sur certaines in\'egalit\'es d'interpolation et sur un r\'esultat r\'ecent de rigidit\'e pour des \'equations elliptiques non lin\'eaires sur ces vari\'et\'es.
\end{abstract}
\end{frontmatter}\vspace*{-1.5cm}
\selectlanguage{english}

\section{Spectral properties of Schr\"odinger operators on the sphere}\label{Sec:sphere}

We start by briefly reviewing some results that have been obtained in \cite{DoEsLa2012}. Let us define $2^*:=\frac{2\,d}{d-2}$ if $d\ge3$, and $2^*:=\infty$ if $d=1$ or $2$.  We denote by $\Lap$ the Laplace-Beltrami operator on  the unit sphere $\S^d\subset\R^{d+1}$. It is well known (see \cite{MR1134481}) that the equation
\[
-\Delta_g u+\tfrac\lambda{q-2}\,u=u^{q-1}
\]
has only constant solutions as long as $q\in(2,2^*)$ and $\lambda\le d$. See \cite{DEKL} for a review and various related results. Assume that the measure on $\S^d$ is the one induced by Lebesgue's measure on $\R^{d+1}$. This convention differs from the one of \cite{DoEsLa2012}. The inequality
\[\label{InterpSphere}
\nrm{\nabla u}2^2+\alpha\,\nrm u2^2\ge\mu(\alpha)\,\nrm uq^2\quad\forall\,u\in\H^1(\S^d)\,,
\]
for any $q\in(2,2^*)$ can be established by standard variational methods. According to \cite{DoEsLa2012}, the optimal function $\mu:\R^+\to\R^+$ is concave, increasing,  and  such that $\mu(\alpha)=\kappaS\,\alpha$ for any $\alpha\le\tfrac d{q-2}$, $\mu(\alpha)<\kappaS\,\alpha$ for $\alpha>\tfrac d{q-2}$ where $\kappaS:=|\S^d|^{1-2/q}$ is a normalization factor and
\[
\mu(\alpha)\sim\mathsf K_{q,d}\,\alpha^{1-\vartheta}\quad\mbox{as}\quad\alpha\to+\infty\,,\quad\mbox{where}\quad\vartheta:=d\,\frac{q-2}{2\,q}\,,
\]
\be{Kqd}
\mathsf K_{q,d}:=\inf_{v\in\H^1(\R^d)\setminus\{0\}}\frac{\nrmRd{\nabla v}2^2+\nrmRd v2^2}{\nrmRd vq^2}\;.
\ee
Let us define $p=\tfrac q{q-2}$ so that $p\in(1,+\infty)$ if $d=1$ and $p\in(\frac d2,+\infty)$ if $d\ge2$. If we denote by $\mu\mapsto\alpha(\mu)$ the inverse function of $\alpha\mapsto\mu(\alpha)$ and by $\lambda_1(-\Delta_g-V)$ the lowest (nonpositive) eigenvalue of $-\Delta_g-V$, then we have the estimate
\[\label{T11}
|\lambda_1(-\Delta_g-V)|\le\alpha\big(\nrm Vp\big)\quad\forall\,V\in\L^p(\S^d)\,.
\]
for any nonnegative $V\in\L^p(\S^d)$. Moreover we have $\alpha(\mu)^{p-d/2}=\L_{p-\frac d2,d}^1\,\mu^p\,(1+o(1))$ as~$\mu\to+\infty$ where $\L_{\gamma,d}^1:=\(\mathsf K_{q,d}\)^{-(\gamma+d/2)}$. Equality is achieved for any $\mu>0$ by some nonnegative $V$, which is constant if and only if $\mu\le\frac d2\,(p-1)$.

The case $q\in(1,2)$ can also be covered and we refer to \cite{DoEsLa2012} for further details. This case leads to estimates from below for the first eigenvalue of the operator $-\Delta_g+W$, where $W$ is a positive potential.

\section{A rigidity result on compact manifolds and a subcritical interpolation inequality}\label{Sec:Rigidity}

{}From here on we shall assume that $(\M,g)$ is a smooth compact connected Riemannian manifold of dimension $d\ge1$, without boundary, and let $\Lap$ be the Laplace-Beltrami operator on $\M$. We shall denote by $d\kern1pt v_g$ the volume element and by~$\Ric$ the Ricci tensor. Let~$\lambda_1$ be the lowest positive eigenvalue of $-\Lap$. 

For such manifolds a new rigidity result has been recently established in \cite{DEL2013}, thus extending a series of results of  \cite{MR615628,MR1134481,MR1338283,MR1412446,MR1631581}. In order to state this result let us define the quantities:
\[
\theta=\frac{(d-1)^2\,(p-1)}{d\,(d+2)+p-1}\quad\mbox{and}\quad\mathrm Q_g u:=\mathrm H_gu-\frac gd\,\Delta_gu-\frac{(d-1)\,(p-1)}{\theta\,(d+3-p)}\left[\frac{\nabla u\otimes\nabla u}u-\frac gd\,\frac{|\nabla u|^2}u\right]
\]
where $\mathrm H_gu$ denotes Hessian of $u$, and
\be{LambdaStar}
\Lambda_\star:=\inf_{u\in\H^2(\M)\setminus\{0\}}\frac{(1-\theta)\iM{(\Delta_gu)^2}+\frac{\theta\,d}{d-1}\iM{\Big[\,\|\mathrm Q_g u\|^2+\Ric(\nabla u,\nabla u)\Big]}}{\iM{|\nabla u|^2}}\,.
\ee
It is not difficult to see that $\Lambda_\star\le\lambda_1$.
\smallskip\begin{thm}{\rm \cite[\emph{cf.} Theorem 3]{DEL2013}}\label{rigiditythm} 
Assume that  $\Lambda_\star$ is strictly positive. Then for any $q\in(1,2)\cup(2,2^*)$ and
any $\lambda\in(0,\Lambda_\star)$, the equation
\[\label{Eqn}
-\,\Lap v+\frac\lambda{q-2}\,\(v-v^{q-1}\)=0
\]
has $1$ as its unique positive solution in $C^2(\M)$. \end{thm}\smallskip
Note that in the particular case $\M=\S^d$, $\Lambda_\star=\lambda_1(-\Lap)=d$. The proof relies on the nonlinear flow
\be{flow}
u_t=u^{2-2\,\beta} \(\Delta_gu+\big(1+\beta\,(q-2)\big)\,\frac{|\nabla u|^2}u\)\,,\quad\beta=\frac{(d+2)\,(d+3-p)\,\theta}{(d-1)^2\,(p-1)^2-(d+2)^2\,(p-2)\,\theta}\;,
\ee
that can also be used to prove the following  $A$--$B$ type interpolation inequality (see \cite{MR1134481,MR1688256}). Let us define
\[\label{Kappa}
\kappaM:=\mathrm{vol}_g(\M)^{1-2/q}\,.
\]

\smallskip\begin{thm}{\rm \cite[\emph{cf.} Theorem 4]{DEL2013}}\label{Thm:Inequality} For any $q\in(1,2)\cup(2,2^*]$ if $d\ge3$, $q\in(1,2)\cup(2,\infty)$ if $d=1$ or $2$, the inequality
\[\label{Ineq:Interp}
\nrmM{\nabla v}2^2\ge\frac\Lambda{q-2}\,\left[\kappaM\,\nrmM vq^2-\nrmM v2^2\right]\quad\forall\,v\in\H^1(\M)\,.
\]
holds for some optimal $\Lambda\in[\Lambda_\star,\lambda_1]$ if $\Lambda_\star>0$. Moreover, if $\Lambda_\star<\lambda_1$, then we have $\Lambda_\star<\Lambda\le\lambda_1$.\end{thm}\smallskip
The above results hold true because the flow~\eqref{flow} contracts
\[\label{functional}
\mathcal F[u]:=\iM{|\nabla(u^\beta)|^2}+\frac{\Lambda_\star}{q-2}\left[\iM{u^{2\,\beta}}-\kappa\,{\(\iM{u^{\beta\,q}}\)}^{2/q}\right]\,.
\]
The above choices for $\theta$ and $\beta$ are optimal for this contraction property: see \cite{DEL2013}.\medskip

As a consequence and exactly as in the case of the sphere, we get  the first  result of this note.
\smallskip\begin{prop}\label{Prop:Interpolation} Assume that $q\in(2,2^*)$ if $d\ge 3$, or $q\in(2,\infty)$ if $d=1$ or $2$. There exists a concave increasing function $\muM:\R^+\to\R^+$ such that $\mu(\alpha)=\kappaM\,\alpha$ for any $\alpha\le\tfrac\Lambda{q-2}$, $\mu(\alpha)<\kappaM\,\alpha$ for $\alpha>\tfrac\Lambda{q-2}$~and
\[\label{InterpM}
\nrmM{\nabla u}2^2+\alpha\,\nrmM u2^2\ge\muM(\alpha)\,\nrmM uq^2\quad\forall\,u\in\H^1(\M)\,.
\]
The asymptotic behaviour of $\muM$ is given by $\muM(\alpha)\sim\mathsf K_{q,d}\,\alpha^{1-\vartheta}$ as $\alpha\to+\infty$, with $\vartheta=d\,\frac{q-2}{2\,q}$ and $\mathsf K_{q,d}$ defined by \eqref{Kqd}.
\end{prop}\smallskip
\begin{proof} There is an optimal function for the interpolation inequality, as can be shown by standard variational techniques. Applying Theorem~\ref{rigiditythm} to the solutions of the Euler-Lagrange equations completes the proof for fixed values of $\alpha$. As an infimum on $u$ of affine functions with respect to $\alpha$,  the function $\alpha\to\mu(\alpha)$ is  concave. It remains to establish the properties of $\alpha$ for large values of $\mu$.

Using a well chosen test function obtained by scaling an optimal function for \eqref{Kqd} on the tangent plane to an arbitrary point of $\M$, one can prove that $\limsup_{\alpha\to+\infty}\alpha^{\vartheta-1}\,\muM(\alpha)\le\mathsf K_{q,d}$. Arguing by contradiction as in \cite[Proposition 10]{DoEsLa2012}, we can find a sequence $\seq\alpha n$ such that $\lim_{n\to+\infty}\alpha_n=+\infty$ and $
\lim_{n\to+\infty}\alpha_n^{\vartheta-1}\mu(\alpha_n)<\mathsf K_{q,d}$, and a sequence of optimal functions $\seq un$ such that $\nrm{u_n}q=1$, which concentrates because $\limsup_{n\to+\infty}\alpha_n^\vartheta\nrm{u_n}2^2<\mathsf K_{q,d}$. Some classical surgery and a convexity inequality provide a contradiction by constructing a minimizing sequence for \eqref{Kqd}.\end{proof}

\section{Ground state estimates for Schr\"odinger operators on Riemannian manifolds}

In this section, we keep using the notations of Section~\ref{Sec:Rigidity} and generalize to $(\M,g)$ the spectral results established for the sphere in \cite{DoEsLa2012}. By inverting the function $\alpha\mapsto \muM(\alpha)$, we see that $\alphaM:\R_+\to\R_+$ is increasing, convex and satisfies: $\alphaM(\mu)=\frac\mu\kappaM$ for any $\mu\in\big(0,\frac{\kappa\,\Lambda}{q-2})$, $\alphaM(\mu)>\frac\mu\kappaM$ for $\mu\in(\frac{\kappa\,\Lambda}{q-2},+\infty)$. With $\L_{\gamma,d}^1:=\(\mathsf K_{q,d}\)^{-p}$, $\gamma=p-\frac d2$, we obtain for a general manifold $\M$ the same behavior of $\mu\mapsto\alphaM(\mu)$ when $\mu\to+\infty$ as in the case of a sphere.

Let $\mu:=\nrm Vp$. Since $p$ and $\frac q2$ are H\"older conjugate exponents, it follows from H\"older's inequality that
\[
\iM{|\nabla u|^2}-\iM{V\,|u|^2}+\alphaM(\mu)\,\iM{|u|^2}\ge\nrm{\nabla u}2^2-\mu\,\nrm uq^2+\alphaM(\mu)\,\nrm u2^2
\]
with equality if $V^{p-1}$ and $|u|^2$ are proportional. The right-hand side is nonnegative according to Proposition~\ref{Prop:Interpolation}. By taking the infimum of the left-hand side, we can deduce an estimate of the lowest, nonpositive eigenvalue $\lambda_1(-\Delta_g-V)$ of $-\Delta_g-V$, which provides us with our first main result.
\smallskip\begin{thm}\label{Thm2} Let $d\ge1$, $p\in(1,+\infty)$ if $d=1$ and $p\in(\frac d2,+\infty)$ if $d\ge2$ and assume that $\Lambda_\star>0$. With the above notations and definitions, for any nonnegative $V\in\L^p(\M)$, we have
\be{T11M}
|\lambda_1(-\Delta_g-V)|\le\alphaM\big(\nrmM Vp\big)\,.
\ee
Moreover, we have $\alphaM(\mu)^{p-\frac d2}=\L_{p-\frac d2,d}^1\,\mu^p\,(1+o(1))$ as $\mu\to+\infty$.\end{thm}
The estimate \eqref{T11M} is optimal in the sense that for any $\mu\in(0,+\infty)$, there exists a nonnegative function~$V$ such that $\mu= \nrmM Vp$ and $|\lambda_1(-\Delta_g-V)|=\alphaM(\mu)$ . Moreover, if $\mu<\frac{\kappa\,\Lambda\star}{q-2}$, $\alphaM(\mu)=\frac\mu\kappaM$ and equality in~\eqref{T11M} is achieved by constant potentials.

\medskip In the case of operators $-\Lap+W$ on $\M$, where $W$ is a nonnegative potential, following again the same arguments as in \cite{DoEsLa2012} in the case of the sphere, and the rigidity result of Theorem~\ref{rigiditythm}, we obtain our second main result.
\smallskip\begin{thm}\label{Thm4} Let $d\ge1$, $p\in(0,+\infty)$. There exists an increasing concave function $\nu:\R^+\to\R^+$, satisfying $\nu(\beta)=\beta/\kappa$, for any $\beta\in(0,\frac{p+1}2\,\kappa\,\Lambda)$ if $p>1$, such that for any positive potential $W$ we have
\[\label{T31}
\lambda_1(-\Delta+W)\ge\nu\big(\beta\big)\quad\mbox{with}\quad\beta={\textstyle\(\iM{W^{-p}}\)^{1/p}}\,.
\]
Moreover, for large values of $\beta$, we have $\nu(\beta)^{-\,(p+\frac d2)}=\L_{-(p+\frac d2),d}^1\,\beta^{-p}\,(1+o(1))$ as $\beta\to+\infty$.\end{thm}

With $p=\tfrac q{2-q}$, the spectral estimate of Theorem~\ref{Thm4} is derived from the interpolation inequality
\[\label{Interpmanifold2}
\nrmM{\nabla u}2^2+\beta\,{\textstyle\(\iM{|u|^q}\)^{2/q}}\ge\nu(\beta)\,\nrmM u2^2\quad\forall\,u\in\H^1(\M)\,.
\]
The concentration phenomena leading to the asymptotics for large norms of $W$  can be studied as in Proposition~\ref{Prop:Interpolation}: see \cite{DoEsLa2012} for the proof in the case of a sphere.
We omit the details of the proof of Theorem~\ref{Thm4}.

\begin{spacing}{0.8}
\linespread{0.9}

\end{spacing}

\smallskip{\small
\noindent{\bf Acknowledgments.} J.D.~and M.J.E.~have been partially supported by ANR grants \emph{CBDif} and \emph{NoNAP}. They thank the Mittag-Leffler Institute, where part of this research was carried out, for hospitality. M.L.~was supported in part by NSF grant DMS-0901304.\\[4pt]
{\sl\small \copyright~2013 by the authors. This paper may be reproduced, in its entirety, for non-commercial purposes.}
\begin{flushright}{\sl\today}\end{flushright}
}

\begin{thebibliography}{1}

\bibitem{MR1412446}
{\sc D.~Bakry and M.~Ledoux}, {\em Sobolev inequalities and {M}yers's diameter
  theorem for an abstract {M}arkov generator}, Duke Math. J., 85 (1996),
  pp.~253--270.

\bibitem{MR1134481}
{\sc M.-F. Bidaut-V{\'e}ron and L.~V{\'e}ron}, {\em Nonlinear elliptic
  equations on compact {R}iemannian manifolds and asymptotics of {E}mden
  equations}, Invent. Math., 106 (1991), pp.~489--539.

\bibitem{DEKL}
{\sc J.~Dolbeault, M.~J. Esteban, M.~Kowalczyk, and M.~Loss}, {\em Sharp
  {I}nterpolation {I}nequalities on the {S}phere: {N}ew {M}ethods and
  {C}onsequences}, Chin. Ann. Math. Ser. B, 34 (2013), pp.~99--112.

\bibitem{DoEsLa2012}
{\sc J.~Dolbeault, M.~J. Esteban, and A.~Laptev}, {\em Spectral estimates on
  the sphere}.
\newblock Preprint hal-00770755, to appear in Analysis and PDE.

\bibitem{DEL2013}
{\sc J.~Dolbeault, M.~J. Esteban, and M.~Loss}, {\em Nonlinear flows and
  rigidity results on compact manifolds}.
\newblock Preprint hal-00784887.

\bibitem{MR615628}
{\sc B.~Gidas and J.~Spruck}, {\em Global and local behavior of positive
  solutions of nonlinear elliptic equations}, Comm. Pure Appl. Math., 34
  (1981), pp.~525--598.

\bibitem{MR1688256}
{\sc E.~Hebey}, {\em Nonlinear analysis on manifolds: {S}obolev spaces and
  inequalities}, vol.~5 of Courant Lecture Notes in Mathematics, New York
  University Courant Institute of Mathematical Sciences, New York, 1999.

\bibitem{MR1338283}
{\sc J.~R. Licois and L.~V{\'e}ron}, {\em Un th\'eor\`eme d'annulation pour des
  \'equations elliptiques non lin\'eaires sur des vari\'et\'es riemanniennes
  compactes}, C. R. Acad. Sci. Paris S\'er. I Math., 320 (1995),
  pp.~1337--1342.

\bibitem{MR1631581}
{\sc J.~R. Licois and L.~V{\'e}ron}, {\em A class of nonlinear conservative
  elliptic equations in cylinders}, Ann. Scuola Norm. Sup. Pisa Cl. Sci. (4),
  26 (1998), pp.~249--283.

\end{thebibliography}
\end{document}